\pdfoutput=1
\documentclass[preprint,1p]{elsarticle}
\usepackage{graphics}
\usepackage{epsfig}
\usepackage{mathptmx}
\usepackage{amsmath}
\usepackage{amssymb}
\usepackage{hyperref}
\usepackage{fullpage}
\usepackage{amsthm}
\usepackage{lineno}
\journal{}
\begin{document}
\begin{frontmatter}
\title{Approximation properties of $(p,q)$-variant of Stancu-Schurer operators(Revised)}
\author[]{Abdul Wafi}
\ead{awafi@jmi.ac.in}
\author[]{Nadeem Rao\corref{cor1}}
\ead{nadeemrao1990@gmail.com}
\address{Department of Mathematics, Jamia Millia Islamia, New Delhi-110 025, India}
\cortext[cor1]{Corresponding author}
\begin{abstract}
In this article, we have introduced $(p,q)$-variant of Stancu-Schurer operators and discussed the rate of convergence for continuous functions. We have also discussed recursive estimates, Korovkin and direct approximation results using second order modulus of continuity, Peetre's K-functional, Lipschitz class.
\end{abstract}
\end{frontmatter}
\section{Introduction}
In 1885, Weierstrass gave a very famous result known as Weierstrass approximation theorem which plays an important role in the development of approximation theory. It was considered to be typical until Bernstein gave an elegant  proof of it. Bernstein\cite{Bernstein} considered polynomials for the continuous functions $f\in C[0,1]$ defined as follows
\begin{equation*}
B_n(f;x)=\sum\limits_{k=0}^{n}P_{n,k}(x)f\Big(\frac{k}{n}\Big), \hspace{2 cm} k=0,1,2,...,n=1,2,3,...
\end{equation*}
where, $P_{n,k}(x)={n\choose k}x^k(1-x)^{n-k}$ and $x\in[0,1]$. It is a  powerful tool for numerical analysis, computer added geometric design(CAGD) and solutions of differential equations.\\

For the last two decades, the application of $q$-calculus emerged as a new area in the field of approximation theory. Motivated by the applications of $q$-calculus, Lupas\cite{Lupas} introduced a sequence of Bernstein polynomials based on $q$-integer. Another form of $q$-Bernstein operators were given by Philips\cite{Philips}. Several researchers introduced different type of operators based on $q$-integers (\cite{ali1}-\cite{gupta1}). Recently, Mursaleen et al and Acar applied $(p,q)$-calculus in approximation theory and introduced $(p,q)$-analogue of Bernstein operators\cite{MB}, Bernstein-Kantorovich operators\cite{MBK}, Bernstein-Stancu operators\cite{MST} and Szasz-Mirakjan operators\cite{tuncer} respectively. The aim of $(p,q)$-integers was to generalize several forms of $q$-oscillator algebras in the earlier physics litrature\cite{chakar}.\\

Let $0<q<p\leq 1$. Then, $(p,q)$-integers for non negative integers $n,k$ are given by
\begin{eqnarray*}
[k]_{p,q}=\frac{p^k-q^k}{p-q}\hspace{1 cm} and\hspace{0.5 cm} [k]_{p,q}=1\hspace{1 cm} for\hspace{0.5 cm} k=0.
\end{eqnarray*}
 $(p,q)$-binomial coefficient
\begin{eqnarray*}
{n \choose k}_{p,q}=\frac{[n]_{p,q}!}{[k]_{p,q}![n-k]_{p,q}!},
\end{eqnarray*}
and $(p,q)$-binomial expansion
\begin{eqnarray*}
(ax+by)_{p,q}^{n}&=&\sum\limits_{k=0}^{n} {n \choose k}_{p,q}p^{\frac{(n-k)(n-k-1)}{2}}q^{\frac{k(k-1)}{2}}a^{n-k}b^kx^{n-k}y^k,\\
(x+y)_{p,q}^{n}&=&(x+y)(px+py)(p^2x+q^2y)...(p^{n-1}x-q^{n-1}y).
\end{eqnarray*}
Mursaleen et al (\cite{MSH}) defined Bernstein-Schurer operators in the following way
\begin{eqnarray}
 B_{n,l}^{p,q}(f;x)=\sum\limits_{\nu=0}^{n+l}b_{n,l}^{\nu}(x;p,q)f\Bigg(\frac{[\nu]_{p,q}}{p^{\nu-(n+l)}[n]_{p,q}}\Bigg) \hspace{2 cm} l=0,1,2,3,...,n=1,2,3,...
\end{eqnarray}
 where $b_{n,l}^{\nu}(x;p,q)=\frac{1}{p^{\frac{(n+l)(n+l-1)}{2}}}{n+l\choose \nu}_{p,q}p^{\frac{\nu(\nu-1)}{2}}x^{\nu}\prod\limits_{j=0}^{n+l-\nu-1}(p^j-q^jx)$.\\
Motivated by the above developments, we define $(p,q)$-variant of Stancu-Schurer operators for any \\
$f\in C[0,l+1], x\in[0,1]$ and $0\leq\alpha\leq\beta$, as follows\\
\begin{eqnarray}
 S_{n,l}^{\alpha,\beta}(f;x,p,q)=\sum\limits_{\nu=0}^{n+l}b_{n,l}^{\nu}(x;p,q)f\Bigg(\frac{p^{(n+l)-\nu}[\nu]_{p,q}+\alpha}{[n]_{p,q}+\beta}\Bigg) \hspace{2 cm} \nu=0,1,2,3,...,n=1,2,3,...
\end{eqnarray}
One can notice that\\
(i) for $\alpha=\beta=0$, (2) reduce to (1).\\
(ii) for $\alpha=\beta=0$ and $p=1$, (2) and reduces to $q$-Stancu-Schurer\cite{pn} operators.
\newline
In the present paper, we investigate the rate of convergence and Korovkin type theorem the operators defined by (2). Moreover, we discussed direct approximation results using second order modulus of continuity, Peeter's K-functional and Lipschitz class of functions.
\label{}
\label{sec2}\section{Basic estimates for Schurer-Stancu operators  $S_{n,l}^{\alpha,\beta}(f;x,p,q)$}
\textbf{Lemma 2.1}\cite{MSH} Let $B_{n,l}^{p,q}(f;x)$ be given by (1). Then for any $x\in[0,1]$ and $0<q<p\leq 1$ we have the following identities
\begin{eqnarray*}
 B_{n,l}^{p,q}(1;x)&=&1,\\
 B_{n,l}^{p,q}(t;x)&=&\frac{[n+l]_{p,q}x}{[n]_{p,q}},\\
 B_{n,l}^{p,q}(t^{2};x)&=&\frac{[n+l]_{p,q}p^{n+l-1}x}{[n]_{p,q}^2}+\frac{q[n+l]_{p,q}[n+l-1]_{p,q}x^2}{[n]_{p,q}^2}.
\end{eqnarray*}
\textbf{Lemma 2.2} Let $x\in[0,1]$ and $0<q<p\leq1$. For the operators $S_{n,l}^{\alpha,\beta}(f;x,p,q)$, we have
\begin{eqnarray*}
 S_{n,l}^{\alpha,\beta}(t^{m};x,p,q)=\frac{[n]_{p,q}^{m}}{([n]_{p,q}+\beta)^{m}}\sum\limits_{i=0}^{m}{m \choose i}\Bigg(\frac{\alpha}{[n]_{p,q}}\Bigg)^{m-i}B_{n,l}^{p,q}(t^{i};x).
\end{eqnarray*}
\textbf{Proof} From (2), we get
\begin{eqnarray*}
 S_{n,l}^{\alpha,\beta}(t^m;x,p,q)&=&\sum\limits_{\nu=0}^{n+l}b_{n,l}^{\nu}(x;p,q)\Bigg(\frac{p^{(n+l)-\nu}[\nu]_{p,q}+\alpha}{[n]_{p,q}+\beta}\Bigg)^{m}\\
 &=&\sum\limits_{\nu=0}^{n+l}b_{n,l}^{\nu}(x;p,q)\frac{[n]_{p,q}^{m}}{([n]_{p,q}+\beta)^{m}}\Bigg(\frac{p^{(n+l)-\nu}[\nu]_{p,q}+\alpha}{[n]_{p,q}}\Bigg)^{m}\\
 &=&\frac{[n]_{p,q}^{m}}{([n]_{p,q}+\beta)^{m}}\sum\limits_{\nu=0}^{n+l}b_{n,l}^{\nu}(x;p,q)\sum\limits_{i=0}^{m}{m \choose i}\Bigg(\frac{\alpha}{[n]_{p,q}}\Bigg)^{m-i}\Bigg(\frac{p^{(n+l)-\nu}[\nu]_{p,q}}{[n]_{p,q}}\Bigg)^{i}\\
 &=&\frac{[n]_{p,q}^{m}}{([n]_{p,q}+\beta)^{m}}\sum\limits_{i=0}^{m}{m \choose i}\Bigg(\frac{\alpha}{[n]_{p,q}}\Bigg)^{m-i}\sum\limits_{\nu=0}^{n+l}b_{n,l}^{\nu}(x;p,q)\Bigg(\frac{p^{(n+l)-\nu}[\nu]_{p,q}}{[n]_{p,q}}\Bigg)^{i}\\
 S_{n,l}^{\alpha,\beta}(t^{m};x,p,q)&=&\frac{[n]_{p,q}^{m}}{([n]_{p,q}+\beta)^{m}}\sum\limits_{i=0}^{m}{m \choose i}\Bigg(\frac{\alpha}{[n]_{p,q}}\Bigg)^{m-i}B_{n,l}^{p,q}(t^{i};x).
\end{eqnarray*}
\textbf{Lemma 2.3} For $S_{n,l}^{\alpha,\beta}$, we have
\begin{eqnarray*}
 S_{n,l}^{\alpha,\beta}(1;x,p,q)&=&1,\\
 S_{n,l}^{\alpha,\beta}(t;x,p,q)&=&\frac{[n+l]_{p,q}x+\alpha}{[n]_{p,q}+\beta},\\
 S_{n,l}^{\alpha,\beta}(t^{2};x,p,q)&=&\frac{[n+l]_{p,q}p^{n+l-1}+2\alpha)x}{([n]_{p,q}+\beta)^2}+\frac{q[n+l]_{p,q}[n+l-1]_{p,q}x^2}{([n]_{p,q}+\beta)^2}+\frac{\alpha^2}{([n]_{p,q}+\beta)^2}.
\end{eqnarray*}
\textbf{Proof} We can prove these equalities using Lemma 2.1 and Lemma 2.2.\\ \\
\textbf{Lemma 2.4} Let $\psi_x^i(t)=(t-x)^i$ and $S_{n,l}^{\alpha,\beta}$ be the operators defined by (3). Then, we get
\begin{eqnarray*}
S_{n,l}^{\alpha,\beta}(\psi_x^0(t);x,p,q)&=&1,\\
S_{n,l}^{\alpha,\beta}(\psi_x^1(t);x,p,q)&=&\Bigg(\frac{[n+l]_{p,q}}{[n]_{p,q}+\beta}-1\Bigg)x+\frac{\alpha}{[n]_{p,q}+\beta},\\
S_{n,l}^{\alpha,\beta}(\psi_x^2(t);x,p,q)&=&\frac{[n+l]_{p,q}[n+l-1]_{p,q}q-2[n+l]_{p,q}([n]_{p,q}+\beta)+([n]_{p,q}+\beta)^2}{([n]_{p,q}+\beta)^2}x^2\\
&&+\frac{[n+l]_{p,q}(p^{n+l-1}+2\alpha)-2\alpha([n]_{p,q}+\beta)}{([n]_{p,q}+\beta)^2}x+\frac{\alpha^2}{(n+\beta)^2}.
\end{eqnarray*}
\section{Convergence properties of $S_{n,l}^{\alpha,\beta}(f;x,p_n,q_n)$}
 \textbf{Theorem 3.1} Let $f\in C[0,l+1]$, $(q_n)_n, (p_n)_n$ be two sequences such that $0<q_n<p_n\leq1$ $ and\lim\limits_{n\rightarrow\infty}{p_n}=1$, \\
 $\lim\limits_{n\rightarrow\infty}{q_n}=1$. Then
 \begin{eqnarray*}
 \parallel S_{n,l}^{\alpha,\beta}(f;x,p_n,q_n)-f(x)\parallel_{C[0,l+1]}=0.
 \end{eqnarray*}
 \textbf{Proof} It is sufficient to show that for $i=0,1,2 $
 \begin{eqnarray*}
  \parallel S_{n,l}^{\alpha,\beta}(f;x,p_n,q_n)(t^i;x)-x^i\parallel\rightarrow 0 \hspace{1 cm} as \hspace{1 cm} n\rightarrow\infty.
 \end{eqnarray*}
 For $i=0$, it is obvious. \\
 For $i=1$, we have
 \begin{eqnarray*}
  |S_{n,l}^{\alpha,\beta}(t;x,p_n,q_n)-x|&=&\Bigg|\frac{[n+l]_{p_n,q_n}x+\alpha}{[n]_{p_n,q_n}+\beta}-x\Bigg|\\
  &\leq& \Bigg(\frac{[n+l]_{p_n,q_n}}{[n]_{p_n,q_n}}-1\Bigg)x+\frac{\alpha}{[n]_{p_n,q_n}+\beta}\\
 \parallel S_{n,l}^{\alpha,\beta}(t;x,p_n,q_n)-x\parallel|&\rightarrow 0& \hspace{1 cm}as \hspace{0.5 cm} n\rightarrow\infty.\\
\end{eqnarray*}
For $i=2$,
 \begin{eqnarray*}
  |S_{n,l}^{\alpha,\beta}(t^2;x,p_n,q_n)-x^2|&=&\Bigg|\frac{[n+l]_{p_n,q_n}(p_n^{n+l-1}+2\alpha)x}{([n]_{p_n,q_n}+\beta)^2}+\frac{q_n[n+l]_{p_n,q_n}[n+l-1]_{p_n,q_n}x^2}{([n]_{p_n,q_n}+\beta)^2}+\frac{\alpha^2}{([n]_{p_n,q_n}+\beta)^2}-x^2\Bigg|\\
  &&\leq\Bigg|\frac{[n+l]_{p_n,q_n}(p_n^{n+l-1}+2\alpha)x}{([n]_{p_n,q_n}+\beta)^2}\Bigg|+\Bigg|\frac{q_n[n+l]_{p_n,q_n}[n+l-1]_{p_n,q_n}}{([n]_{p_n,q_n}+\beta)^2}-1\Bigg|x^2\\
  &&+\Bigg|\frac{\alpha^2}{([n]_{p_n,q_n}+\beta)^2}\Bigg|\\
  ||S_{n,l}^{\alpha,\beta}(t^2;x,p_n,q_n)-x^2||&\rightarrow 0& \hspace{1 cm}as \hspace{0.5 cm}n \rightarrow \infty.
  \end{eqnarray*}
and this proves the theorem. \\ \\
\textbf{Theorem 3.2} Let $(p_n)_n,(q_n)_n$ be the sequences such that $0<q_n<p_n\leq 1$ and $\lim\limits_{n\rightarrow\infty}{p_n}=1$, $\lim\limits_{n\rightarrow\infty}{q_n}=1$. Then
\begin{eqnarray*}
|S_{n,l}^{\alpha,\beta}(f;x,p_n,q_n)-f(x)|\leq 2\omega(f;\sqrt{\delta_{n,l}^{\alpha,\beta}(x)}),
\end{eqnarray*}
for all $f\in C[0,l+1]$ and $\delta_{n,l}^{\alpha,\beta}(x)=\sqrt{S_{n,l}^{\alpha,\beta}(\psi_x^2(t);x,p_n,q_n)}$.\\ \\
\textbf{Proof} Calculating the difference, we find
\begin{eqnarray*}
|S_{n,l}^{\alpha,\beta}(f;x,p_n,q_n)-f(x)|&=&\sum\limits_{\nu=0}^{n+l}b_{n,l}^{\nu}(x;p,q)\Bigg|f\Bigg(\frac{p_n^{(n+l)-\nu}[\nu]_{p_n,q_n}+\alpha}{[n]_{p_n,q_n}+\beta}\Bigg)-f(x)\Bigg|\\
&\leq&\Bigg\{1+\frac{1}{\delta_n^{\alpha,\beta}}\sum\limits_{\nu=0}^{n+l}b_{n,l}^{\nu}(x;p,q)\Bigg|\frac{p_n^{(n+l)-\nu}[\nu]_{p_n,q_n}+\alpha}{[n]_{p_n,q_n}+\beta}-x\Bigg|\Bigg\}\omega(f;\delta_n^{\alpha,\beta}(x))\\
&\leq&\Bigg\{1+\frac{1}{\delta_n^{\alpha,\beta}}\sqrt{\sum\limits_{\nu=0}^{n+l}b_{n,l}^{\nu}(x;p,q)\Big(\bigg|\frac{p_n^{(n+l)-\nu}[\nu]_{p_n,q_n}+\alpha}{[n]_{p_n,q_n}+\beta}-x\bigg|\Big)^2}\Bigg\}\omega(f;\delta_n^{\alpha,\beta}(x))\\
&=&\Bigg\{1+\frac{1}{\delta_n^{\alpha,\beta}}\sqrt{S_{n,l}^{\alpha,\beta}(\psi_x^2(t);x,p_n,q_n)}\Bigg\}\omega(f;\delta_n^{\alpha,\beta}(x))
\end{eqnarray*}
choosing $\delta_{n,l}^{\alpha,\beta}=\sqrt{S_{n,l}^{\alpha,\beta}(\psi_x^2(t);x,p_n,q_n)}$, we get the desired result.\\ \\
 \textbf{Theorem 3.3} For $0\leq\alpha\leq\beta$,  $0<q_n<p_n\leq 1$ such that $\lim\limits_{n\rightarrow\infty}q_n=1$, $\lim\limits_{n\rightarrow\infty}q_n=1,$ and $ x\in [0,1],$ we have
\begin{eqnarray*}
 | S_{n,l}^{\alpha,\beta}(f;x)-f(x) | \leq\omega_1(([n]_{p,q}+\beta)^{-1})\sqrt{S_{n,l}^{\alpha,\beta}(\psi_x^2(t);x)}\bigg\{1+\sqrt{([n]_{p,q}+\beta)}\sqrt{S_{n,l}^{\alpha,\beta}(\psi_x^2(t);x)} \bigg\}
\end{eqnarray*}
where $f'(x)$ has continuous derivative over $[0,l+1]$ and $\omega_1(f;\delta_{n,\beta})$ is the modulus of continuity of $f'(x)$.\\ \\ \\
\textbf{Proof} For $t_1,t_2\in [0,b]$ and $t_1<\eta<t_2$, it is known that
\begin{eqnarray}
\nonumber f(t_1)-f(t_2)&=&(t_1-t_2)f'(\eta),\\
 &=&(t_1-t_2)f'(t_1)+(t_1-t_2)[f'(\eta)-f'(t_1)],
\end{eqnarray}
 and from Lorentz[\cite{lorentz}, p. 21, Theorem 1.6.2], we have
\begin{eqnarray}
 |(t_1-t_2)[f'(\eta)-f'(t_1)]| \leq| t_1-t_2 |(\lambda+1)\omega_1(\delta),
\end{eqnarray}
where $\lambda=\lambda(x_1,x_2;\delta)$ is the integer $[|t_1-t_2|\delta^{-1}]$.\\
Now, we have
\begin{eqnarray}
| S_{n,l}^{\alpha,\beta}(f;x)-f(x) |=\bigg|\sum\limits_{\nu=0}^{n+l}b_{n,l}^{\nu}(x;p,q)f\Bigg(\frac{p_n^{(n+l)-\nu}[\nu]_{p_n,q_n}+\alpha}{[n]_{p_n,q_n}+\beta}\Bigg)-f(x)\bigg|.
\end{eqnarray}
Using (3) and (4), we get
\begin{eqnarray*}
| S_{n,l}^{\alpha,\beta}(f;x)-f(x)| &\leq &\Bigg| \sum\limits_{\nu=0}^{n+l}b_{n,l}^{\nu}(x;p,q)\Bigg(x-\frac{p_n^{(n+l)-\nu}[\nu]_{p_n,q_n}+\alpha}{[n]_{p_n,q_n}+\beta}\Bigg)f'(x)\Bigg|\\
&&+\omega_1(\delta_{n,l}^{\alpha,\beta})\sum\limits_{\nu=0}^{n+l}\Bigg|\frac{p_n^{(n+l)-\nu}[\nu]_{p_n,q_n}+\alpha}{[n]_{p_n,q_n}+\beta}-x\Bigg|(\lambda+1)b_{n,l}^{\nu}(x;p,q),\\ &\leq & \omega_1(\delta_{n,l}^{\alpha,\beta})\Bigg\{\sum\limits_{\nu=0}^{n+l}\Bigg|\frac{p^{(n+l)-\nu}[\nu]_{p_n,q_n}+\alpha}{[n]_{p_n,q_n}+\beta}-x\Bigg|b_{n,l}^{\nu}(x;p,q)\\
&&+\sum\limits_{\nu=0}^{n+l}\Bigg|\frac{p^{(n+l)-\nu}[\nu]_{p_n,q_n}+\alpha}{[n]_{p_n,q_n}+\beta}-x\Bigg|\Bigg(\lambda\Big(x,\frac{[\nu]_{p_n,q_n}+\alpha}{[n]_{p_n,q_n}+\beta};\delta_{n,l}^{\alpha,\beta}\Big)+1\Bigg)b_{n,l}^{\nu}(x;p,q)\Bigg\},\\
&\leq & \omega_1(\delta_{n,l}^{\alpha,\beta})\Bigg\{\sum\limits_{\nu=0}^{n+l}\Bigg|\frac{p^{(n+l)-\nu}[\nu]_{p_n,q_n}+\alpha}{[n]_{p_n,q_n}+\beta}-x\Bigg|b_{n,l}^{\nu}(x;p,q)\\
&&+\frac{1}{\delta_n^\beta}\sum\limits_{\nu=0}^{n+l}\Bigg(\frac{p^{(n+l)-\nu}[\nu]_{p_n,q_n}+\alpha}{[n]_{p_n,q_n}+\beta}-x\Bigg)^2b_{n,l}^{\nu}(x;p,q)\\
&\leq&\omega_1(\delta_n^\beta)\Bigg(\sqrt{S_{n,l}^{\alpha,\beta}(\psi_x^2;x)}+\frac{S_{n,l}^{\alpha,\beta}(\psi_x^2;x)}{\delta_n^\beta}\Bigg)\\
&=&\omega_1(\delta_n^\beta)\sqrt{S_{n,l}^{\alpha,\beta}(\psi_x^2;x)}\Bigg\{1+\frac{\sqrt{S_{n,l}^{\alpha,\beta}(\psi_x^2;x)}}{\delta_n^\beta}\Bigg\}
\end{eqnarray*}
Taking $\delta_n^\beta=([n]_{p,q}+\beta)^{-1}$, we get
\begin{eqnarray*}
 | S_{n,l}^{\alpha,\beta}(f;x)-f(x) | \leq\omega_1(([n]_{p,q}+\beta)^{-1})\sqrt{S_{n,l}^{\alpha,\beta}(\psi_x^2(t);x)}\bigg\{1+\sqrt{([n]_{p,q}+\beta)}\sqrt{S_{n,l}^{\alpha,\beta}(\psi_x^2(t);x)} \bigg\}.
\end{eqnarray*}
\section{Direct results for $S_{n,l}^{\alpha,\beta}(f;x)$ }

Let $C_B[0,\infty)$ denote the space of real valued continuous and bounded functions $f$ on $[0,\infty)$ endowed with the norm
\begin{eqnarray*}
\|f\|=\sup\limits_{0\leq x<\infty}|f(x)|.
\end{eqnarray*}
Then, for any $\delta>0$, Peeter's K-functional is defined as
\begin{eqnarray*}
K_2(f,\delta)=inf\{\|f-g\|+\delta\|g''\|: g\in C_B^2[0,\infty)\},
\end{eqnarray*}
where $C_B^2[0,\infty)=\{g\in C_B[0,\infty):g',g''\in C_B[0,\infty)\}$. Now, we know that there exits an absolute constant $C>0$ Devore and Lorentz[\cite{devor}, p.177, Theorem 2.4] such that
\begin{eqnarray*}
K_2(f;\delta)\leq C\omega_2(f;\sqrt{\delta}),
\end{eqnarray*}
where $\omega_2(f;\delta)$ is the second order modulus of continuity given by
\begin{eqnarray*}
\omega_2(f,\sqrt{\delta})=\sup\limits_{0<h\leq\sqrt{\delta}} \sup\limits_{x\in[0,\infty)} |f(x+2h)-2f(x+h)+f(x)|.
\end{eqnarray*}
\textbf{Theorem 4.1} Let $f\in C_B^2[0,l+1]$, $(q_n)_n, (p_n)_n$ be two sequences such that $0<q_n<p_n\leq1$ $ and\lim\limits_{n\rightarrow\infty}{p_n}=1$. Then for all $x\in[0,\infty)$ there exists a constant $K>0$ such that
\begin{eqnarray*}
\mid S_{n,l}^{\alpha,\beta}(f;x)-f(x)\mid\leq K\omega_2(f;\sqrt{\Theta_{n,a}^{\alpha,\beta}(x)})+\omega\Bigg(f;\Bigg(\frac{[n+l]_{p_n,q_n}}{[n]_{p_n,q_n}+\beta}-1\Bigg)x+\frac{\alpha}{[n]_{p_n,q_n}+\beta}\Bigg)
\end{eqnarray*}
where
\begin{eqnarray*}
\Theta_{n,a}^{\alpha,\beta}(x)&=&S_{n,l}^{\alpha,\beta}((t-x)^2;x)+\Bigg(\bigg(\frac{[n+l]_{p_n,q_n}}{[n]_{p_n,q_n}+\beta}-1\bigg)x+\frac{\alpha}{[n]_{p_n,q_n}+\beta}\Bigg)^2.\\
\end{eqnarray*}
\textbf{Proof} First, we consider the auxiliary operators $ \widehat{S}_{n,l}^{\alpha,\beta}$
\begin{eqnarray}
\widehat{S}_{n,l}^{\alpha,\beta}(f;x)=S_{n,l}^{\alpha,\beta}(f;x)+f(x)-f\Bigg(\bigg(\frac{[n+l]_{p_n,q_n}}{[n]_{p_n,q_n}+\beta}-1\bigg)x+\frac{\alpha}{[n]_{p_n,q_n}+\beta}\Bigg).
\end{eqnarray}
We find that
\begin{eqnarray*}
\widehat{S}_{n,l}^{\alpha,\beta}(1;x)=1,
\end{eqnarray*}
\begin{eqnarray*}
\widehat{S}_{n,l}^{\alpha,\beta}(t-x;x)=0,
\end{eqnarray*}
\begin{eqnarray}
|\widehat{S}_{n,l}^{\alpha,\beta}(f;x)|\leq 3\|f\|.
\end{eqnarray}
Let $g\in C_B^2[0,\infty)$. By the Taylor's theorem\\
\begin{eqnarray}
g(t)=g(x)+(t-x)g'(x)+\int\limits_x^t (t-v)g''(v)dv.
\end{eqnarray}
Now
\begin{eqnarray*}
\widehat{S}_{n,a}^{\alpha,\beta}(g;x)-g(x)&=&g'(x)\widehat{S}_{n,a}^{\alpha,\beta}(t-x;x)+\widehat{S}_{n,l}^{\alpha,\beta}\Big( \int\limits_x^t (t-v)g''(v)dv;x\Big)\\
&=&\widehat{S}_{n,a}^{\alpha,\beta}\Big( \int\limits_x^t (t-v)g''(v)dv;x\Big)\\
&=&S_{n,a}^{\alpha,\beta}\Big( \int\limits_x^t (t-v)g''(v)dv;x\Big)\\
&&-\int\limits_x^{\big(\frac{[n+l]_{p_n,q_n}}{[n]_{p_n,q_n}+\beta}-1\big)x+\frac{\alpha}{[n]_{p_n,q_n}+\beta}}\Bigg(\bigg(\frac{[n+l]_{p_n,q_n}}{[n]_{p_n,q_n}+\beta}-1\bigg)x+\frac{\alpha}{[n]_{p_n,q_n}+\beta}-v\Bigg)g''(v)dv.
\end{eqnarray*}
Therefore
\begin{eqnarray}
\nonumber | \widehat{S}_n^{\alpha,\beta}(g;x)-g(x)|&\leq&\Bigg|S_{n,a}^{\alpha,\beta}\Big( \int\limits_x^t (t-v)g''(v)dv;x\Big)\Bigg|\\
&&+\Bigg|\int\limits_x^{\big(\frac{[n+l]_{p_n,q_n}}{[n]_{p_n,q_n}+\beta}-1\big)x+\frac{\alpha}{[n]_{p_n,q_n}+\beta}} \Bigg(\bigg(\frac{[n+l]_{p_n,q_n}}{[n]_{p_n,q_n}+\beta}-1\bigg)x+\frac{\alpha}{[n]_{p_n,q_n}+\beta}-v\Bigg)g''(v)dv\Bigg|.
\end{eqnarray}
Since
 \begin{eqnarray}
\Bigg| \int\limits_x^t (t-v)g''(v)dv\Bigg|\leq(t-x)^2\parallel g''\parallel,
 \end{eqnarray}
 then we have
 \begin{eqnarray}
 \nonumber \Bigg|\int\limits_x^{\big(\frac{[n+l]_{p_n,q_n}}{[n]_{p_n,q_n}+\beta}-1\big)x+\frac{\alpha}{[n]_{p_n,q_n}+\beta}} \Big(\bigg(\frac{[n+l]_{p_n,q_n}}{[n]_{p_n,q_n}+\beta}-1\bigg)x+\frac{\alpha}{[n]_{p_n,q_n}+\beta}-v\Big)g''(v)dv\Bigg|&\leq&\Big(\bigg(\frac{[n+l]_{p_n,q_n}}{[n]_{p_n,q_n}+\beta}-1\bigg)x\\
 &&+\frac{\alpha}{[n]_{p_n,q_n}+\beta}\Big)^2\parallel g''\parallel.
 \end{eqnarray}
 From(6),(7) and (8), we have
 \begin{eqnarray}
\nonumber | \widehat{S}_{n,l}^{\alpha,\beta}(g;x)-g(x)|&\leq& \Bigg\{S_{n,l}^{\alpha,\beta}((t-x)^2;x)+\bigg(\frac{[n+l]_{p_n,q_n}}{[n]_{p_n,q_n}+\beta}-1\bigg)x+\frac{\alpha}{[n]_{p_n,q_n}+\beta}\bigg)^2\Bigg\}\|g''\|\\
 &=&\Pi_{n,l}^{\alpha,\beta}(x)\|g''\|.
 \end{eqnarray}
 Next, we have
 \begin{eqnarray*}
 |S_{n,l}^{\alpha,\beta}(f;x)-f(x)|&\leq& |\widehat{S}_{n,l}^{\alpha,\beta}(f-g;x)|+|(f-g)(x)|+|\widehat{S}_{n,l}^{\alpha,\beta}(g;x)-g(x)|\\
 &&+\Big|f\Bigg(\bigg(\frac{[n+l]_{p_n,q_n}}{[n]_{p_n,q_n}+\beta}-1\bigg)x+\frac{\alpha}{[n]_{p_n,q_n}+\beta}\bigg)\Bigg)-f(x)\Big|
 \end{eqnarray*}
 Using(9), we have
 \begin{eqnarray*}
 |S_{n,l}^{\alpha,\beta}(f;x)-f(x)|&\leq& 4\|f-g\| +\widehat{S}_{n,l}^{\alpha,\beta}(g;x)-g(x)|+\bigg|f\Bigg(\bigg(\frac{[n+l]_{p_n,q_n}}{[n]_{p_n,q_n}+\beta}-1\bigg)x+\frac{\alpha}{[n]_{p_n,q_n}+\beta}\bigg)-f(x)\bigg|\\
 &\leq& 4\|f-g\|+\Theta_{n,l}^{\alpha,\beta}(x)\|g''\|+\omega\Bigg(f;\bigg(\frac{[n+l]_{p_n,q_n}}{[n]_{p_n,q_n}+\beta}-1\bigg)x+\frac{\alpha}{[n]_{p_n,q_n}+\beta}\Bigg).
 \end{eqnarray*}
 By the definition of Peetre's K-functional, we have
 \begin{eqnarray*}
 |S_{n,l}^{\alpha,\beta}(f;x)-f(x)|\leq C\omega_2\big(f;\sqrt{\Theta_{n,l}^{\alpha,\beta}(x)}\big)+\omega(f;\bigg(\frac{[n+l]_{p_n,q_n}}{[n]_{p_n,q_n}+\beta}-1\bigg)x+\frac{\alpha}{[n]_{p_n,q_n}+\beta}\Bigg).
 \end{eqnarray*}
  Now, we discuss the rate of convergence of the operators $S_{n,l}^{\alpha,\beta}$ in Lipschitz class $Lip_M (\alpha)$, given by
 \begin{eqnarray*}
 Lip_M(\alpha)=\{f\in C[0,\infty):|f(t)-f(x)|\leq M|t-x|^\alpha:x,t\in[0,1]\}
 \end{eqnarray*}
 where $M$ is a constant and $0<\alpha\leq 1$.\\ \\ \\
 \textbf{Theorem 4.3} Let $f\in Lip_M(\alpha)$ and $0<q<p\leq 1$. Then, we have
 \begin{eqnarray*}
  |S_{n,l}^{\alpha,\beta}(f;x)-f(x)|\leq M\delta_{n,l}^{\alpha,\beta}(x),
 \end{eqnarray*}
 where $\delta_{n,l}^{\alpha,\beta}(x)=S_{n,l}^{\alpha,\beta}((t-x)^2;x)$. \\ \\
 \textbf{Proof} Let $\alpha=1$ and $0<q<p\leq 1$. Then for $f\in Lip_M(1)$, we have
 \begin{eqnarray*}
  |S_{n,l}^{\alpha,\beta}(f;x)-f(x)|&\leq & \sum\limits_{\nu=0}^{n+l}b_{n,l}^{\nu}(x;p,q)\Bigg|f\Bigg(\frac{p_n^{(n+l)-\nu}[\nu]_{p,q}+\alpha}{[n]_{p,q}+\beta}\Bigg)-f(x)\Bigg| \\
  &\leq & M\sum\limits_{\nu=0}^{n+l}b_{n,l}^{\nu}(x;p,q)\Bigg|\frac{p_n^{(n+l)-\nu}[\nu]_{p,q}+\alpha}{[n]_{p,q}+\beta}-x\Bigg| \\
  &= & M S_{n,l}^{\alpha,\beta}(|t-x|;x)\\
  &\leq& M (S_{n,l}^{\alpha,\beta}((t-x)^2;x))^{\frac{1}{2}}\\
  &=&M (\delta_{n,l}^{\alpha,\beta}(x))^{\frac{1}{2}}
 \end{eqnarray*}
 Thus, the assertion hold for $\alpha=1$. Now, we will prove for $\alpha\in (0,1)$. From the Holder inequality with $p=\frac{1}{\alpha}, \\ q=\frac{1}{1-\alpha}$, we have
 \begin{eqnarray*}
  |S_{n,l}^{\alpha,\beta}(f;x)-f(x)|&=&\Bigg( \sum\limits_{\nu=0}^{n+l}b_{n,l}^{\nu}(x;p,q)\Bigg|f\Bigg(\frac{p_n^{(n+l)-\nu}[\nu]_{p,q}+\alpha}{[n]_{p,q}+\beta}\Bigg)-f(x)\Bigg|^{\alpha}\Bigg)^{\frac{1}{\alpha}}\\
  &&\times\Bigg(\sum\limits_{\nu=0}^{n+l}b_{n,l}^{\nu}(x;p,q)\Bigg)^{\frac{1}{1-\alpha}}\\
  &\leq& \Bigg( \sum\limits_{\nu=0}^{n+l}b_{n,l}^{\nu}(x;p,q)\Bigg|f\Bigg(\frac{p_n^{(n+l)-\nu}[\nu]_{p,q}+\alpha}{[n]_{p,q}+\beta}\Bigg)-f(x)\Bigg|^{\alpha}\Bigg)^{\frac{1}{\alpha}}
 \end{eqnarray*}
 Since $f\in Lip_M(\alpha)$, we obtain
 \begin{eqnarray*}
  |S_{n,l}^{\alpha,\beta}(f;x)-f(x)|&\leq& M\Bigg(\sum\limits_{\nu=0}^{n+l}b_{n,l}^{\nu}(x;p,q)\Bigg|\frac{p_n^{(n+l)-\nu}[\nu]_{p,q}+\alpha}{[n]_{p,q}+\beta}-x\Bigg|\Bigg)^\alpha\\ &\leq& M(S_{n,l}^{\alpha,\beta}(|t-x|;x))^{\alpha}\\
  &\leq& M (S_{n,l}^{\alpha,\beta}((t-x)^2;x))^{\frac{\alpha}{2}}\\
  &=&M (\delta_{n,l}^{\alpha,\beta}(x))^{\frac{\alpha}{2}}
  \end{eqnarray*}

 \textbf{Remark} Approximation results obtained for Bernstein-Schurer by Mursaleen at al(see\cite{MSH}) are the particular case of our results for $\alpha=\beta=0$.

\end{document}